\documentclass[onefignum,onetabnum]{siamart171218}

\usepackage{amsfonts}
\usepackage{graphicx}
\usepackage{epstopdf}
\usepackage{algorithmic}
\ifpdf
  \DeclareGraphicsExtensions{.eps,.pdf,.png,.jpg}
\else
  \DeclareGraphicsExtensions{.eps}
\fi

\newsiamremark{remark}{Remark}
\newsiamremark{hypothesis}{Hypothesis}
\crefname{hypothesis}{Hypothesis}{Hypotheses}
\newsiamthm{claim}{Claim}

\headers{Scalable higher-order nonlinear solvers via higher-order AD}{S. Tan, K. Miao, A. Edelman, and C. Rackauckas}

\title{Scalable higher-order nonlinear solvers via higher-order automatic differentiation}

\author{
  Songchen Tan\thanks{Department of Mathematics, Massachusetts Institute of Technology, Cambridge, MA 
  (\email{songchen@mit.edu}).}
  \and
  Keming Miao
  \and
  Alan Edelman
  \and
  Christopher Rackauckas\thanks{Department of Mathematics, Massachusetts Institute of Technology, Cambridge MA 
  (\email{crackauc@mit.edu}).}
}

\usepackage{amsopn}

\makeatletter
\newcommand*{\addFileDependency}[1]{
  \typeout{(#1)}
  \@addtofilelist{#1}
  \IfFileExists{#1}{}{\typeout{No file #1.}}
}
\makeatother

\ifpdf
\hypersetup{
  pdftitle={An Example Article},
  pdfauthor={D. Doe, P. T. Frank, and J. E. Smith}
}
\fi

\begin{document}

\maketitle

\begin{abstract}
  This paper demonstrates new methods and implementations of nonlinear solvers with higher-order of convergence, which is achieved by efficiently computing higher-order derivatives. Instead of computing full derivatives, which could be expensive, we compute directional derivatives with Taylor-mode automatic differentiation. We first implement Householder's method with arbitrary order for one variable, and investigate the trade-off between computational cost and convergence order. We find that the second-order variant, i.e., Halley's method, to be the most valuable, and further generalize Halley's method to systems of nonlinear equations and demonstrate that it can scale efficiently to large-scale problems. We further apply Halley's method on solving large-scale ill-conditioned nonlinear problems, as well as solving nonlinear equations inside stiff ODE solvers, and demonstrate that it could outperform Newton's method.
\end{abstract}



\section{Introduction}

Solving nonlinear equations accurately and efficiently is crucial in many disciplines of science, where nonlinear equations can emerge either directly from the models, or as an intermediate step of solving ordinary differential equations (ODEs)\cite{byrne1987stiff}, partial differential equations (PDEs)\cite{ames2014numerical}, differential-algebraic equations (DAEs)\cite{kunkel2006differential}, and integral equations\cite{atkinson1992survey}.

Various methods are available to solve nonlinear equations, and the mostly used ones are Newton's method and its variants\cite{kelley2003solving}, which rely on the first-order derivative, i.e. the Jacobian, or its approximations. Newton's method is proved to have quadratic convergence under certain assumptions of the nonlinear function and initial condition\cite{galantai2000theory}.

There has been extensive research on improving the convergence order of nonlinear solvers, which can be categorized into two different approaches: (1) multi-step methods that evaluate nonlinear functions and/or derivatives at multiple points in each iteration\cite{potra1984nondiscrete,xiao_new_2015,madhu2017some,xiao_new_2022,singh_simple_2023}; (2) methods that utilize higher-order derivatives of nonlinear functions\cite{alefeld1981convergence,cuyt_abstract_1983,cuyt_computational_1985,cuyt1982numerical,noor2007modified,ogbereyivwe2023some,sariman2020new,germani_higher-order_2006}. Despite their higher convergence order, in order to be more efficient than Newton's method, they must be computationally efficient for each iteration. For example, some multi-step methods\cite{singh_simple_2023} evaluate the function and Jacobian at multiple points, but only invert the Jacobian at one point; in practice, this would mean that for a nonlinear function with $n$ inputs and $n$ outputs, the $O(n^3)$ factorization of the Jacobian, which is the bottleneck of Newton's method, only needs to be done once for each iteration. This ensures that these multi-step methods are not much more expensive than Newton's method for each iteration.

The same considerations also apply to methods of the second category, where higher-order derivatives introduce additional computation burden. It is commonly believed that computing the higher-order full derivative of a function with large input and output dimension $n$ would be much more expensive than the computation of the function itself\cite{margossian2019review}. As a result, they are often approximated rather than exactly computed\cite{sariman2020new,suleiman2009solving,albeanu2008generalized,steihaug2006newton}, and it is still an open question whether nonlinear solvers can utilize higher-order derivative information for large-scale nonlinear problems in an exact way.

This paper demonstrates new methods and implementations that, instead of computing the full derivative, efficiently and exactly compute the higher-order directional derivative, which is sufficient for nonlinear solvers. This paper is organized as follows. We will first introduce the technique of Taylor-mode automatic differentiation (AD) for efficiently computing higher-order directional derivatives, and its implementation in Julia, in \cref{sec:theory}. Then we will demonstrate that Taylor-mode AD can be used to efficiently implement Householder's method with arbitrary convergence order, in \cref{sec:householder}. Finally, we demonstrate that a special case of Householder's method, the cubic-convergence Halley's method can scale elegantly to large-scale problems, and could outperform first-order methods by a significant margin, in \cref{sec:experiments}.

\section{Taylor-mode automatic differentiation for higher-order directional derivatives}
\label{sec:theory}

We first introduce the notation used in this paper. Let $U$ be an open subset of $\mathbb R^n$, $f:U\to\mathbb R^m$ be a function that is sufficiently smooth. The derivative of $f$ at a point $x\in\mathbb R^n$ is a linear operator $Df(x):\mathbb R^n\to\mathbb R^m$ such that it maps a vector $v\in\mathbb R^n$ to $Df(x)[v]$, which is the directional derivative of $f$ at $x$ in the direction of $v$, also known as the Jacobian-vector product.

Similarly, for $p\in\mathbb Z^+>1$, the $p$-th order derivative of $f$ at a point $x$ is a multilinear operator $D^pf(x):(\mathbb R^n,\cdots,\mathbb R^n)\to\mathbb R^m$ such that it maps a tuple of vectors $(v_1,\ldots,v_p)$ to $D^pf(x)[v_1,\cdots,v_p]$, the directional derivative of $f$ at $x$ in the direction of $(v_1,\ldots,v_p)$.

First-order forward-mode AD can be viewed as an algorithm to propagate directional derivative information through compositions of functions\cite{revels_forward-mode_2016}. For example, let $f(x)=g(h(x))$ to be a composite function for $\mathbb R^n\to\mathbb R^m$, and $x\in\mathbb R^n$ to be a specific point in its domain. Suppose that we already have $h_0=h(x)$ to represent the primal output of $h$, and $h_1=Dh(x)[v]$ representing the perturbation of $h$ along some direction $v\in\mathbb R^n$. Now for the composite function, we want to know what is the perturbation of $f$ along direction $v$; in other words, we want to know $f_1=Df(x)[v]$ in addition to $f_0=f(x)$. The answer to that is simply the chain rule:
\begin{itemize}
  \item $f_0 = g(h_0)$
  \item $f_1 = Dg(h_0)[h_1]$
\end{itemize}
Similarly, the Taylor-mode AD can be viewed as an algorithm to propagate higher-order directional derivative information through compositions of functions. Suppose that we already have a Taylor bundle for $h$ at $x$ along $v$:
$$
(h_0,h_1,\ldots,h_p) = (h(x),Dh(x)[v],D^2h(x)[v,v],\ldots,D^ph(x)[v,\ldots,v])
$$
and we want to know
$$
(f_0,f_1,\ldots,f_p) = (f(x),Df(x)[v],D^2f(x)[v,v],\ldots,D^pf(x)[v,\ldots,v])
$$
the answer to that is, in turn, the Fa\`a di Bruno's formula:
\begin{itemize}
  \item $f_0 = g(h_0)$
  \item $f_1 = Dg(h_0)[h_1]$
  \item $f_2 = D^2g(h_0)[h_1,h_1] + Dg(h_0)[h_2]$
  \item $\cdots$
\end{itemize}
In our practical implementation of Taylor-mode AD, a pushforward rule is defined for every ``simple'' function, so that derivatives of any complicated functions that are composed of these simple functions can be computed automatically, either via operator-overloading or source-code transformation techniques. Our previous work \cite{tan_higher-order_2023} has shown that for any functions composed of elementary functions and arbitrary control flow, the $p$-th order directional derivative is only $O(p)$ times more expensive than computing the function itself, in contrast to naively nesting first-order forward-mode AD which could be $O(\exp(p))$ times expensive\cite{bettencourt_taylor-mode_2022}. We also provided a reference implementation of this method in the Julia language\cite{bezanson_julia_2017} using its multiple-dispatch mechanism\cite{bezanson_fast_2018}, TaylorDiff.jl\cite{tan_taylordiffjl_2023}, featuring an automatic generation\cite{gowda_high-performance_2021} of higher-order pushforward rules from first-order pushforward rules defined in ChainRules.jl\cite{white_juliadiffchainrulesjl_2023}.

Below, we will demonstrate how to use Taylor-mode AD to implement higher-order nonlinear solvers. For all numerical experiments in this paper, we use the Julia language v1.11.2 on a single core of Intel Xeon Platinum 8260 CPU provided by the MIT SuperCloud system.

\section{Efficient implementation of Householder's method}
\label{sec:householder}

We first implement Householder's method\cite{householder1970numerical} for solving nonlinear equations with one variable, i.e. $x\in\mathbb R$ in the nonlinear function $f(x)$. This does not immediately generalize to multiple variables, but it is useful to compare the behavior of solvers with different convergence order. The method is defined as follows: given an initial guess $x_0$, in each iteration we compute,
$$
x_{n+1} = x_n + p\frac{(1/f)^{(p-1)}(x_n)}{(1/f)^{(p)}(x_n)}.
$$
which requires derivatives up to $p$-th order and has convergence order $p+1$. Note that when $p=1$, this method is equivalent to Newton's method; and when $p=2$, this method is equivalent to Halley's method.

Given $f$ and $x_n$, we can compute $f(x_n)$ and its derivatives up to $p$-th order using Taylor mode AD, and then invert this Taylor polynomial to get the value and derivatives of $1/f$ up to $p$-th order. Our analysis leads to the \cref{alg:householder}:

\begin{algorithm}
\caption{Householder's method implementation based on higher-order AD}
\label{alg:householder}
\begin{algorithmic}
\STATE{Define $f$ and initial value $x$; tolerance $t$; order $p$}
\WHILE{$|f(x)|>t$}
\STATE{Initialize a Taylor bundle $(x,1,0,\cdots,0)$}
\STATE{Apply the pushforward rule of $f$ to get $(f(x),f'(x),\cdots,f^{(p)}(x))$}
\STATE{Apply the pushforward rule of $(\cdot)^{-1}$ to get $((1/f)(x),(1/f)'(x),\cdots,(1/f)^{(p)}(x))$}
\STATE{Update $x := x + p\frac{(1/f)^{(p-1)}(x)}{(1/f)^{(p)}(x)}.$}
\ENDWHILE
\RETURN $x$
\end{algorithmic}
\end{algorithm}

\subsection{Cost-effectiveness for Householder's method with different orders.}

\begin{figure}[htbp]
  \centering
  \includegraphics[width=\textwidth]{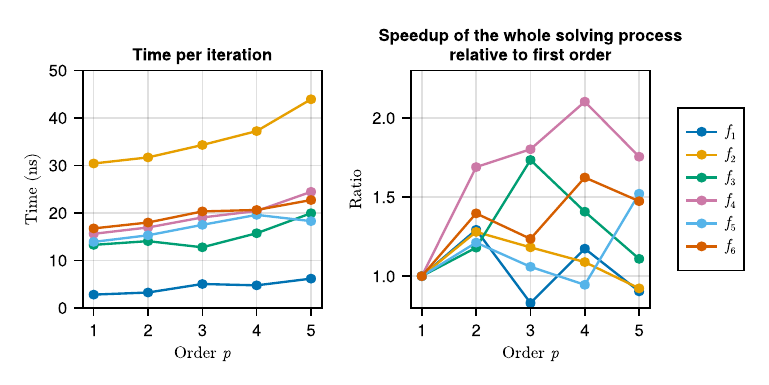}
  \caption{Cost-effectiveness for Householder's method with different orders. Left: the computation time per iteration for each nonlinear function and each Householder's method with different order. Right: the relative speedup for each function of the total computation time to solve the nonlinear equation, normalized by the $p=1$ solver (Newton's method). The functions and initial conditions are: (1) $f_1(x)=x^2-2$, $x_0=1.0$; (2) $f_2(x)=\sqrt{x}-\pi$, $x_0=10.0$; (3) $f_3(x)=x-\exp(-x)$, $x_0=0.0$; (4) $f_4(x)=x^2-2^x$, $x_0=3.3$; (5) $f_5(x)=x+\sin(x)-1$, $x_0=0.5$; (6) $f_6(x)=\log(x)+x$, $x_0=1.0$.}
  \label{fig:householder}
\end{figure}

We implement \cref{alg:householder} as a part of SimpleNonlinearSolve.jl\cite{pal_nonlinearsolvejl_2024}, which is a package of non-allocating nonlinear solvers that have low-overhead and could solve small nonlinear problems very efficiently. We make several univariate nonlinear functions composed of various elementary functions, and use different orders of Householder's method to solve $f(x)=0$ given appropriate initial guess. Results are shown in \cref{fig:householder}.

For each of the nonlinear functions, the computation cost per iteration only slightly increases as the order $p$ goes from 1 (Newton's method) to 5, but the number of iterations needed to converge could be cut down significantly by the higher-order methods. As the result, the total time needed to converge could be shorter for higher-order methods. We also notice that the most significant improvement is from the first order to the second order, and the improvement from the second order to an even higher order is less significant.

\subsection{Generalization to multivariate functions}

The theory of abstract rational approximation\cite{cuyt_abstract_1983} gives a method to generalize Householder's method to multivariate functions. In the sections below, we will focus on implementing $p=2$ case, i.e. Halley's method, with abstract rational approximation formulas, for multivariate functions. Based on the observations from simple univariate problems, we postulate that methods for $p\ge 3$ are less likely to be significantly more efficient than Halley's method.

\section{Scalable and Efficient implementation of Halley's method}
\label{sec:experiments}

Halley's method ($p=2$) for a multivariate real function $f(x):\mathbb R^n\to\mathbb R^n$, derived from abstract rational approximation\cite{cuyt_abstract_1983,cuyt_computational_1985}, could be expressed as follows: given an initial guess $x_0$, in each iteration we compute
\begin{equation}
  x_{n+1}=x_n+(a_n\odot a_n)\oslash (a_n + [Df(x_n)]^{-1}D^2f(x_n)[a_n,a_n] / 2)
\end{equation}
where $\odot$ and $\oslash$ are element-wise multiplication and division, $Df$ and $D^2f$ are the Jacobian and Hessian of $f$, and $a_n$ is the solution to $Df(x_n)a_n=-f(x_n)$. In the equation above, $D^2f(x_n)a_na_n$ can be computed easily using Taylor-mode AD, with just two times more expensive than the computation of $f(x_n)$. After the first linear solve for $a_n$, we need another linear solve to get $[Df(x_n)]^{-1}D^2f(x_n)a_na_n$. Among the two linear solves, the factorization of Jacobian could be reused, similar to the case of multi-step methods\cite{singh_simple_2023}. Our analysis leads to \cref{alg:mvhalley}:

\begin{algorithm}
\caption{(Multivariate) Halley's method implementation based on higher-order AD}
\label{alg:mvhalley}
\begin{algorithmic}
\STATE{Define $f$ and initial value $x$; tolerance $t$}
\WHILE{$\|f(x)\|_{\inf}>t$}
\STATE{Obtain via first-order AD the Jacobian $Df(x)$}
\STATE{Solve $a$ from $Df(x)[a]=-f(x)$}
\STATE{Initialize a Taylor bundle $(x,a,0)$}
\STATE{Apply the pushforward rule of $f$ to get $(f(x),Df(x)[a],D^2f(x)[a,a])$}
\STATE{Solve $b$ from $Df(x)[b]=D^2(x)[a,a]$}
\STATE{Update $x := x + a \odot a \oslash (a + b / 2)$}
\ENDWHILE
\RETURN $x$
\end{algorithmic}
\end{algorithm}

Looking into the details of \cref{alg:mvhalley}, we note that for either dense or sparse Jacobian, the construction and factorization of $Df(x)$ (via LU or QR factorization) should be the bottleneck of solving, and these only need to be done once for the two linear solves. In other words, given the factorization, the additional linear solve for $b$ can be computed in $O(n^2)$ time. Therefore, each step of Halley's method can be asymptotically as cheap as Newton's method, which makes it potentially faster than Newton's method.

\subsection{Solving nonlinear problems with dense Jacobian}

We first test \cref{alg:mvhalley} on nonlinear problems where the Jacobian of the nonlinear function $f(x):\mathbb R^n\to\mathbb R^n$ is dense and unstructured. In this case, the factorization of the Jacobian is the bottleneck of the linear solve, and the cost of each factorization is $O(n^3)$. For defining such a problem, we consider the Chandrasekhar's H-function\cite{chandrasekhar2013radiative} which is the solution to an integral equation that arises in the study of radiative transfer in astrophysics. The equation is defined as follows:
\begin{equation}
  H(\mu)=1+\mu H(\mu)\int_0^1\frac{\Psi(\mu')H(\mu')}{\mu+\mu'}\mathrm d\mu'
\end{equation}
where $H(\mu)$ is defined on $[0,1]$ and $\Psi(\mu)$ is an even polynomial. We consider the most simplified case, where $\Psi$ is a constant. Upon discretization onto $n$ grids for integration, the equation becomes a system of nonlinear equations of $n$ variables, and since each equation depends on the value of $H$ on all positions, the Jacobian of this nonlinear problem is inherently dense. The detailed form of the discretized nonlinear function is available as the 23rd problem in \href{https://github.com/SciML/DiffEqProblemLibrary.jl/blob/master/lib/NonlinearProblemLibrary/src/NonlinearProblemLibrary.jl#L438}{NonlinearProblemLibrary.jl}\cite{pal_nonlinearsolvejl_2024} which the reader can refer to.

\begin{figure}[htbp]
  \centering
  \label{fig:scaling}
  \includegraphics[width=.7\textwidth]{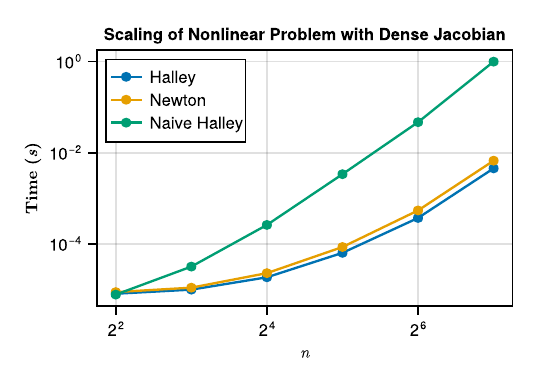}
  \caption{Scaling for solving a nonlinear problem that has a dense Jacobian, at different problem sizes. Solvers: Newton, Halley (Taylor-mode AD), and Naive Halley (full Hessian).}
\end{figure}

We solve this problem with different problem sizes $n=4, 8, 16, 32, 64, 128$, and we compare the performance of (1) Newton's method, (2) Halley's method as implemented in \cref{alg:mvhalley}, and (3) naive Halley's method implemented in a conventional way\cite{cuyt_computational_1985} that computes the full Hessian and then contracts it with vectors. Each of them is solved to the default tolerance. The results are shown in \cref{fig:scaling}.

\begin{figure}[htbp]
  \centering
  \label{fig:wp}
  \includegraphics[width=\textwidth]{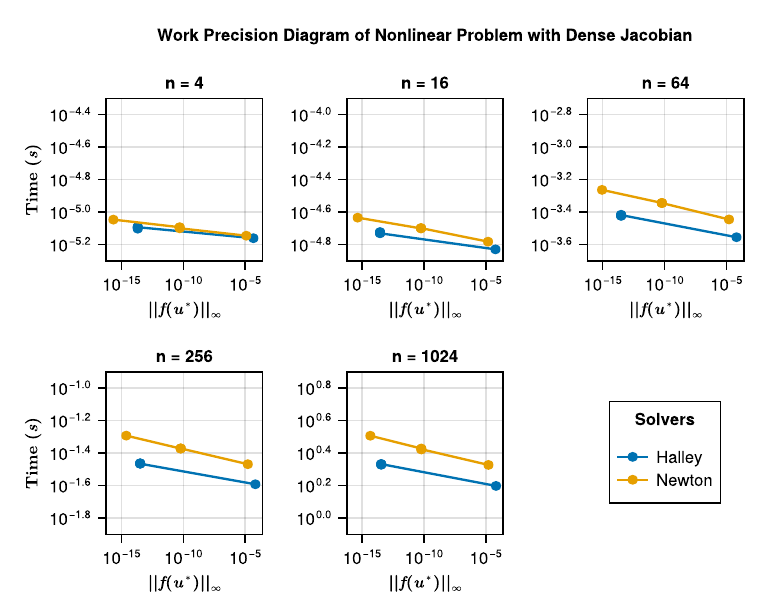}
  \caption{Work-precision diagram for solving a nonlinear problem that has a dense Jacobian, with different sizes. Solvers: Newton and Halley (Taylor-mode AD).}
\end{figure}

While the performances are similar for small problems, the naive Halley's method quickly becomes infeasible for large problems\cite{pal_nonlinearsolvejl_2024}, while our implementation of Halley's method scales similarly as Newton's method. In addition, as the problem size gets larger, the advantage of Halley's method over Newton's method becomes more significant, as a result of fewer iterations needed to converge. For example, for $n=128$, Halley's method is approximately 40\% faster than Newton's method. Finally, more detailed work-precision diagrams comparing Halley's method and Newton's method, within a broader range of problem sizes ($n=4,16,64,256,1024)$, is shown in \cref{fig:wp}.

\subsection{Solving large-scale ill-conditioned nonlinear problems with sparse Jacobian}

Nonlinear problems that appear in solving PDEs are often large, sparse and ill-conditioned\cite{ames2014numerical}, imposing additional challenges for nonlinear solvers. Below, we will use a two-dimensional Brusselator reaction-diffusion PDE problem as an example to demonstrate the capability of Halley's method to handle these nonlinear problems in practical applications. The Brusselator is a model for auto-catalytic chemical systems that exhibit oscillations in time domain and pattern formation in space domain\cite{prigogine1968symmetry}. Let functions $u(x,y,t)$ and $v(x,y,t)$ be concentrations of substances that are defined on $(x, y)\in[0,1]\times[0,1]$ and $t\in\mathbb (0,+\infty)$, satisfying the following differential equations,

\begin{equation}
  \begin{aligned}
    u_t&=B + u^2v - (A+1)u + \alpha\Delta u + f(x,y,t) \\
    v_t&=Au - u^2v + \alpha\Delta v
  \end{aligned}
\end{equation}
where the source $f$ is defined as
\begin{equation}
  f(x,y,t)=\begin{cases}
    5 & \text{if } (x - 0.3)^2 + (y - 0.6)^2 \le 0.1^2 \text{and } t \ge 1.1 \\
    0 & \text{otherwise}
  \end{cases}
\end{equation}
and initial conditions,
\begin{equation}
  \begin{aligned}
    u(x,y,0)&=22(y(1-y))^{3/2} \\
    v(x,y,0)&=27(x(1-x))^{3/2}
  \end{aligned}
\end{equation}
and periodic boundary conditions.
\begin{equation}
  \begin{aligned}
    u(0,y,t)&=u(1,y,t) \\
    u(x,0,t)&=u(x,1,t) \\
    v(0,y,t)&=v(1,y,t) \\
    v(x,0,t)&=v(x,1,t)
  \end{aligned}
\end{equation}

In the following numerical experiments, the parameters are set to $A=3.4$, $B=1$, $\alpha=10$.

We first consider the time-independent PDE problem, where for $t\gg 1.1$ we solve the steady-state equations $u_t=v_t=0$ to get the long-time behavior of the spatial distribution of concentrations $u$ and $v$. In order to apply nonlinear solvers, the spatial discretization of $u$ and $v$ is carried out for $K=4, 8, 16, 32, 64, 128$ grids on $x$ and $y$ directions, and the Laplacian operator $\Delta$ is implemented with finite differences. Therefore, each of $u$ and $v$ discretizes to $K^2$ variables, and the total problem size of the nonlinear system is $n=2K^2$. Since the variables only interact with nearby grids, the resulting Jacobian is sparse; in addition, when $K$ increases, the coefficients of the finite difference operator become larger and make the Jacobian ill-conditioned.

\begin{figure}[htbp]
  \centering
  \label{fig:b}
  \includegraphics[width=\textwidth]{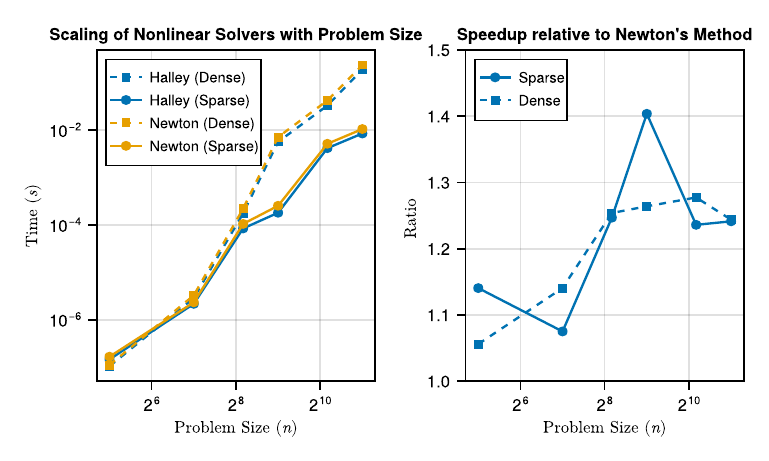}
  \caption{Scaling for solving discretized two-dimensional Brusselator steady-state problem with different problem sizes. The problems are solved to default tolerance. Left: computation time for Newton and Halley's method, each with or without sparsity detection, for each problem size. Right: relative speedup of Halley's method over Newton's method with or without sparsity detection, for each problem size.}
\end{figure}

In order to handle these characteristics, we first try to generate a sparse matrix for Jacobian in order to make factorization easier\cite{averick1994computing}. This could be achieved with sparse detection techniques, either in a numerical way\cite{giering_generating_2005,walther2009getting} or in a symbolic way\cite{gowda_sparsity_2019}. After sparsity detection, we further use the graph coloring algorithm\cite{gebremedhin_what_2005} to compute the matrix. These tools are available in SparseDiffTools.jl\cite{noauthor_juliadiffsparsedifftoolsjl_2024} and SparseConnectivityTracer.jl\cite{hill_sparseconnectivitytracerjl_2024} as a part of the Julia community. Finally, factorization of the sparse Jacobian is appropriately handled with LinearSolve.jl\cite{noauthor_scimllinearsolvejl_2024} which has reasonable polyalgorithms for handling ill-conditioned sparse matrices. In addition, we also include the naive method of computing dense Jacobian and then doing dense factorization. The results are shown in \cref{fig:b}.

Clearly, the naive method of computing dense Jacobian and then doing dense factorization is inefficient for large problems, while the sparse method scales better. In addition, Halley's method is more efficient than Newton's method for all problem sizes and for both sparse and dense handling of the Jacobian matrix, and the advantage becomes more significant as the problem size increases. For example, for $K=32$ and corresponding problem size $n=2048$, Halley's method is approximately 25\% faster than Newton's method.

\subsection{Solving stiff ODEs}

We then consider the time-dependent PDE, where we compute the time evolution of the Brusselator system given initial conditions and periodic boundary conditions above. In order to implement, the PDE is again spatially discretized with $K$ grids on $x$ and $y$ directions, turning it into an ODE. Due to the Laplacian operator, the linear part of this equation has a high frequency, while the nonlinear part has a low frequency. Therefore, the resulting ODE is stiff. In solving stiff ODEs, implicit solvers are often used to ensure accuracy and efficiency\cite{kim_stiff_2021}, and these solvers requires nonlinear solve steps. Since the nonlinear solve steps are often the bottleneck of implicit solvers, impoving the nonlinear solver could accelerate the whole ODE solving process.

We choose problem size $K=8$, $n=2K^2=128$ to construct the ODE and solve with several common stiff solvers from $t_0=0$ to $t_{\rm end}=11.5$, and in each of the stiff solvers we test both Newton's method and Halley's method for the underlying nonlinear solver. When solving the ODE, we provide different tolerances which will determine the time step by applying adaptive time-step strategies available in DifferentialEquations.jl\cite{rackauckas2017differentialequations}. We then use a very small tolerance to obtain very accurate solutions to the ODE, namely $u_{\rm ref}$ and $v_{\rm ref}$, and we define the error to be the relative $L^2$ error of the final solution at $t_{\rm end}$ on domain $D=[0,1]\times[0,1]$, i.e.

\begin{equation}
  \text{error} = \frac{\int_D\left(|u_{\rm ref}(x,y,t_{\rm end})-u(x,y,t_{\rm end})|^2+|v_{\rm ref}(x,y,t_{\rm end})-v(x,y,t_{\rm end})|^2\right)\mathrm dx\mathrm dy}{\int_D\left(|u_{\rm ref}(x,y,t_{\rm end})|^2+|v_{\rm ref}(x,y,t_{\rm end})|^2\right)\mathrm dx\mathrm dy}
\end{equation}

\begin{figure}[htbp]
  \centering
  \label{fig:stiff}
  \includegraphics[width=\textwidth]{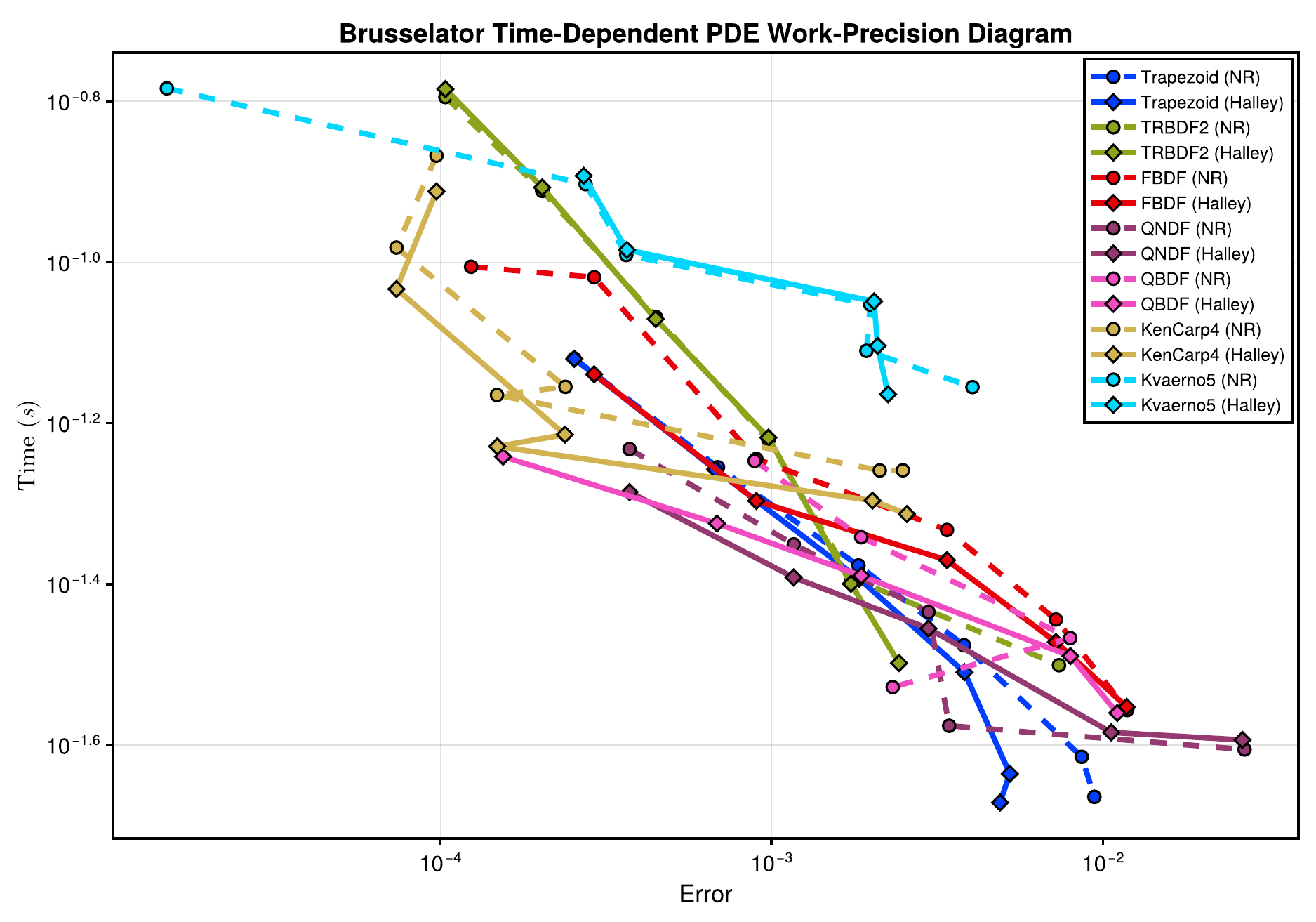}
  \caption{Work-precision diagram for solving discretized two-dimensional Brusselator time-dependent PDE. Implicit Solvers: Trapezoid\cite{vladimirescu_spice_1994}, TRBDF2\cite{hosea_analysis_1996}, FBDF\cite{shampine_solving_2002}, QNDF\cite{shampine_matlab_1997}, QBDF (alias of QNDF with $\kappa=0$), KenCarp4\cite{kennedy_additive_2003}, Kvaerno5\cite{kvaerno_singly_2004}. For each solver, the Newton's method is drawn in circle and solid line, and the Halley's method is drawn in diamond and dashed line.}
\end{figure}

The results are shown in \cref{fig:stiff}.

In the work-precision diagram, we observe that Halley's method is more efficient than Newton's method for many kinds of implicit solvers, and on average it gives a 5-10\% speedup for the whole ODE solving process.

\section{Conclusions}
\label{sec:conclusions}

In this paper, we have demonstrated that higher-order nonlinear solvers can be implemented efficiently and scalable with higher-order automatic differentiation. We have implemented Householder's method with arbitrary convergence order and investigated the trade-off between computational cost and accelerated convergence. We have implemented $p=2$ variant, Halley's method, for multivariate functions and have shown that Halley's method can be more efficient than Newton's method for several cases, including problems with dense Jacobian, problems with sparse Jacobian, as well as being integrated as a part of stiff ODE solvers. Backed by the Julia language and NonlinearSolve.jl framework, this method can be applied to a wide range of problems in science and engineering, and we believe that it could be a valuable tool for the scientific community.

The most noticeable limitation of this work is that we have not yet considered the case where the linear equations in each iteration must be solved iteratively, such as the case that the Jacobian is a matrix-free operator, which is common in even larger systems\cite{KNOLL2004357}. In this case, no factorization could be reused between two linear solves required by Halley's method, and the cost of each iteration could be twice as expensive as Newton's method, cancelling out the advantage of higher convergence order. We will leave this as future work.

\section*{Acknowledgments}
This material is based upon work supported by the U.S. National Science Foundation under award Nos CNS-2346520, PHY-2028125, RISE-2425761, DMS-2325184, OAC-2103804, and OSI-2029670, by the Defense Advanced Research Projects Agency (DARPA) under Agreement No. HR00112490488,  by the Department of Energy, National Nuclear Security Administration under Award Number DE-NA0003965 and by the United States Air Force Research Laboratory under Cooperative Agreement Number FA8750-19-2-1000.  Neither the United States Government nor any agency thereof, nor any of their employees, makes any warranty, express or implied, or assumes any legal liability or responsibility for the accuracy, completeness, or usefulness of any information, apparatus, product, or process disclosed, or represents that its use would not infringe privately owned rights. Reference herein to any specific commercial product, process, or service by trade name, trademark, manufacturer, or otherwise does not necessarily constitute or imply its endorsement, recommendation, or favoring by the United States Government or any agency thereof. The views and opinions of authors expressed herein do not necessarily state or reflect those of the United States Government or any agency thereof. The views and conclusions contained in this document are those of the authors and should not be interpreted as representing the official policies, either expressed or implied, of the United States Air Force or the U.S. Government.
\nocite{*}
\bibliographystyle{siamplain}
\bibliography{references}
\end{document}